\documentclass{amsart}
\usepackage{amssymb}

\newtheorem{theorem}{Theorem}[section]
\newtheorem{lemma}[theorem]{Lemma}
\newtheorem{proposition}[theorem]{Proposition}
\newtheorem{corollary}[theorem]{Corollary}
\newtheorem{question}[theorem]{Question}
\newtheorem{_definition}[theorem]{Definition}
\newenvironment{definition}{\begin{_definition}\rm}{\end{_definition}}
\newtheorem{_remark}[theorem]{\it Remark}
\newenvironment{remark}{\begin{_remark}\rm}{\end{_remark}}
\newtheorem{_claim}[theorem]{Claim}

\numberwithin{equation}{section}
\numberwithin{table}{section}
\numberwithin{figure}{section}
\renewcommand{\P}{\mathord{\mathbb  P}}
\newcommand{\Z}{\mathord{\mathbb Z}}
\newcommand{\R}{\mathord{\mathbb R}}

\newcommand{\C}{\mathord{\mathbb C}}
\newcommand{\Q}{\mathord{\mathbb Q}}
\newcommand{\Pic}{\mathord{\rm{Pic}}}
\newcommand{\Null}{\mathord{\rm{Null}}}
\newcommand{\Stab}{\mathord{\rm{Stab}}}
\newcommand{\Supp}{\mathord{\rm{Supp}}}
\newcommand{\Imm}{\mathord{\rm{Im}}}
\newcommand{\tor}{\mathord{\rm{tor}}}
\newcommand{\NE}{\mathord{\rm{NE}}}
\newcommand{\OO}{\mathord{\mathcal O}}

\newcommand{\Gal}{\mathord{\rm{Gal}}}
\newcommand{\Ker}{\mathord{\rm{Ker}}}

\newtheorem{setup}[theorem]{}

\newcommand{\NS}{\mathord{\rm{NS}}}

\newcommand{\id}{\mathord{\rm{id}}}

\newcommand{\Aut}{\mathord{\rm{Aut}}}

\newcommand{\rank}{\operatorname{\mathord{\rm{rank}}}\nolimits}

\newcommand{\GL}{\mathord{\hbox{\sl{GL}}}}

\begin{document}
\begin{large}

\title
[Automorphism groups]{Automorphism groups and anti-pluricanonical
curves}
%
\author{De-Qi Zhang}
\address
{%
\textsc{Department of Mathematics} \endgraf
\textsc{National University of Singapore, 2 Science Drive 2,
Singapore 117543, Singapore
\endgraf
and
\endgraf
Universit\"at Duisburg-Essen,
Fachbereich Mathematik,
45117 Essen, Germany
}}
\email{matzdq@nus.edu.sg}
%
%
%
%
%
%
%
%
\begin{abstract}
We show the existence of an anti-pluricanonical curve on every
smooth projective rational surface $X$ which has an infinite group $G$ of
automorphisms of either null entropy or of type $\Z \ltimes \Z$,
provided that the pair $(X, G)$ is minimal. This was conjectured by
Curtis T. McMullen (2005) and further traced back to Marat
Gizatullin and Brian Harbourne (1987).
We also prove (perhaps) the strongest form of the famous Tits
alternative theorem.
\end{abstract}
\subjclass{14J26; 14J50, 37B40, 32H50} \keywords{automorphism,
anti-pluricanonical curve, complex dynamics} \maketitle
\section{Introduction}
\par
In this paper, we are interested in the automorphism groups of
smooth projective complex surfaces. Especially, we are interested in the following
question of Gizatullin-Harbourne-McMullen
(see \cite{HPAMS} page 409 and \cite{McP2} \S 12),
where $\Aut^*(X) := \Imm(\Aut(X)$ $\rightarrow \Aut(\Pic(X)))$:
\par \vskip 0.5pc
\begin{question}\label{HQ}
Let $X$ be a smooth projective rational surface.
If $\Aut^*(X)$ is infinite,
is there then a birational morphism $\varphi$ of $X$ to a surface $Y$
having an anti-pluricanonical curve and an infinite subgroup $G \subset \Aut^*(Y)$
such that $G$ lifts via $\varphi$ to $X ?$
\end{question}
\par \vskip 0.5pc
\par
A member in an anti-pluricanonical system
$|-nK_X|$ ($n\ge 1$) is called an {\it anti n-canonical} curve (or {\it divisor})
or simply an {\it anti-pluricanonical} curve;
a member in $|-K_X|$ is an {\it anti 1-canonical} curve,
or an {\it anti-canonical} curve.
\par
The result below answers Question \ref{HQ}
in the case of null entropy.
\par \vskip 0.5pc
\begin{theorem}\label{Thnull}
Let $X$ be a smooth projective rational surface and $G \le \Aut(X)$
an infinite subgroup of null entropy $($see $\ref{entropy})$.
Then we have:
\begin{itemize}
\item[(1)]
There is a $G$-equivariant smooth blowdown $X \rightarrow Y$
such that $K_Y^2 \ge 0$ and hence $Y$ has an anti-pluricanonical curve.
\item[(2)]
Suppose further that $\Imm (G \rightarrow \Aut(\Pic(X)))$
is also an infinite group. Then the $Y$ in $(1)$ can be so chosen that
$-K_Y$ is nef of self intersection zero
and $Y$ has an anti $1$-canonical curve.
\end{itemize}
\end{theorem}
\par \vskip 0.5pc
For groups which are not necessarily of null entropy,
we have the following resul
which is especially applicable (with the same kind of $H$)
when $G/H \ge \Z \ltimes \Z$. See Theorems \ref{Thz2a} - \ref{Thcent}
for more general results.
\par \vskip 0.5pc
\begin{theorem}\label{Thz2}
Let $X$ be a smooth projective surface and $G \le \Aut(X)$ a subgroup.
Assume that there is a sequence of groups
$$H \trianglelefteq A \trianglelefteq G$$
satisfying the following three conditions:
\begin{itemize}
\item[(1)]
$\Imm(H \rightarrow \Aut(\NS(X)))$ is finite;
\item[(2)]
$A/H$ is infinite and abelian; and
\item[(3)]
$|G/A| = \infty$.
\end{itemize}
Then $G$ contains a subgroup $S$ of null entropy and infinite order.
\newline
In particular, when $X$ is rational, there is an $S$-equivariant
smooth blowdown $X \rightarrow Y$ such that $Y$ has an anti-pluricanonical curve.
\end{theorem}
\par \vskip 0.5pc
\begin{remark}
(1) Conditions like the ones in Theorems \ref{Thnull} and \ref{Thz2} are
probably necessary in order to have an affirmative answer to
Question \ref{HQ}. See Bedford-Kim
\cite[Theorem 3.2]{BK} for a pair $(X, g)$
with $g$ of positive entropy and Iitaka D-dimension $\kappa(X, -K_X) = -\infty$.

\par
(2) The blowdown process $X \rightarrow Y$ to the minimal pair $(Y, S)$
in Theorems \ref{Thnull} and \ref{Thz2} is necessary,
as observed by Harbourne \cite{HPAMS}.
\par
However, the minimality assumption of the pair in Theorem \ref{Thpos-a}
is not essential, though it makes the statement simpler.
\end{remark}
\par \vskip 0.5pc
\par
To tackle Question \ref{HQ} of Gizatullin-Harbourne-McMullen,
we begin by determining the dynamical structure of surface automorphism groups.
We first take a close look at McMullen's smooth rational surface $S^v$.
It turns out that the dynamics of $\Aut(S^v)$ is very simple:
\par \vskip 0.5pc
\begin{theorem}\label{Thsv}
Suppose that $v$ is a leading eigenvector for a Coxeter element
$w$ in the Weyl group $W_n \ (n \ge 10)$.
\par
Then for the McMullen surfaces $S^v$ of \cite[\S 7] {McP2} $($the surfaces
$S_n$ in \cite[Theorem 1.1] {McP2} are among them$)$,
we have $($semi-direct product$):$
$$\Aut(S^v) = \langle h_m \rangle \ltimes T$$
with $h_m$ of positive entropy and $T \triangleleft \Aut(S^v)$ a finite subgroup.
\end{theorem}
\par \vskip 0.5pc
Our result below is (perhaps) the strongest form of
the famous Tits \cite{Ti} alternative theorem,
for surface automorphisms groups of positive entropy.
Let $\sigma : \Aut(X) \rightarrow \Aut(\NS_{\Q}(X))$ be the natural
homomorphism.
\par
Note that the $T$ itself in (2) below is also finite when $X$
is a rational or $K3$ surface (see Proposition \ref{cyclo}; Sterk
\cite{St} Lemma 2.1, and Torelli theorem).
\par \vskip 0.5pc
\begin{theorem}\label{Thtits}
Let $X$ be a smooth projective surface and $G \le \Aut(X)$
a subgroup of positive entropy. Then $G$ satisfies either:
\begin{itemize}
\item[(1)]
$G$ contains the non-abelian free group $\Z * \Z$; or
\item[(2)]
There is a $B \trianglelefteq G$ such that
$|G / B| \le 2$ and $B = \langle h_m \rangle \ltimes T$
\newline
$($semi-direct product$)$
with $h_m$ positive entropy and $\sigma(T)$ finite.
\end{itemize}
\end{theorem}
\par \vskip 0.5pc
See Oguiso \cite{Og2} Theorem 2.1 or 1.3
for groups of null entropy and Oguiso \cite{Og1} Theorem 1.1 for $K3$
groups (and more generally for hyperk\"ahler manifolds), where the
cyclic-ness of $B/T$ in our result here is replaced by abelian-ness.
The Case(1) in the theorem above does occur (see Mazur \cite{Ma},
Cantat \cite{Ca3}, Cantat - Favre \cite{CF} Example 3.2, and Oguiso
\cite{Og2} Theorem 1.6).
\par
In Theorem \ref{Thpos-a} below
we determine the relation between
the anti- pluricanonical curve and the set $\Stab(g)$ of $g$-periodic curves.
\par \vskip 0.5pc
\begin{theorem}\label{Thpos-a}
Let $X$ be a smooth projective rational surface with
an automorphism $g$ of positive entropy.
Assume the following two conditions:
\begin{itemize}
\item[(1)] The pair $(X, \langle g \rangle)$ is minimal; and
\item[(2)] Either the set $\Stab(g)$ of $g$-periodic curves contains a curve
of arithmetic genus $\ge 1$,
or $X$ has an anti-pluricanonical curve.
\end{itemize}
Then we have:
\begin{itemize}
\item[(1)]
There is a unique nonzero effective $\Q$-divisor $\Delta$ with
$\Supp(\Delta) \subseteq \Stab(g)$ such that
$K_X + \Delta \equiv 0 \ ($numerical equivalence$)$.
\item[(2)] $\Stab(g)$ is a union of $\Supp(\Delta)$ and possibly
a few $(-2)$-curves away from $\Supp(\Delta)$.
\item[(3)]
Every anti-pluricanonical curve is a multiple of $\Delta$.
\end{itemize}
\end{theorem}
\begin{setup} {\rm New development}
\end{setup}
In \cite{NZ}, \cite{Z3}, \cite{Z4} and \cite{KOZ},
we prove the conjecture of Amerik-Campana-Fujimoto-Nakayama
on the induced action on the base variety of an Iitaka fibration
(generalizing the result of Deligne-Nakamura-Ueno),
prove the conjecture of Guedj on cohomologically hyperbolic endomorphisms, and show
some theorems of Tits type (generalizing results of Dinh-Sibony \cite{DS1}).
\par \vskip 0.5pc
\begin{setup} {\bf Acknowledgment}
\end{setup}
I would like to thank
Florin Ambro for patiently answering my questions on $\R$-divisors,
Serge Cantat and JongHae Keum for constructive comments,
Curtis T. McMullen for informing me about his new paper \cite{McP2},
Noboru Nakayama for answering my questions on
strictly nef divisors, Keiji Oguiso for informing me of his result
in \cite{Og2} used in our Theorem \ref{Thnulla},
Helene Esnault and Eckart Viehweg for the support and warm hospitality
during my stay in Essen,
and the anonymous referee for kind suggestions.
\par
The Main Theorem of \cite{DS2} on automorphisms group $G$
(in the surface case, and with the assumption of the commutativity
of $G$) agrees with the result of our Theorem \ref {Thz2a} on $G$
(which is not necessarily commutative).
\par
This project was supported by an Academic Research Fund of NUS.
\section{Conventions and Preliminary results}
\par \vskip 0.1pc
\begin{setup} {\bf Conventions}
{\rm are as in \cite{H}, \cite{KMM} and \cite{KM}.}
\end{setup}
\par
Let $V$ be a normal projective variety.
For an $\R$-Cartier divisor $D$ on $V$, \,
$\Null(D) := \{C \, | \, C$ is an irreducible curve with $C . D = 0\}$
is the set of curves {\it annihilated} by $D$.
Set
$D^{\perp} := \{E \in \NS_{\Q}(V) \, | \, E . D = 0\}$.

\par
For an automorphism $g \in \Aut(V)$, a curve
$C$ is $g$-{\it periodic} if $g^s(C) = C$ for some $s > 0$.
Denote by $\Stab(g) := \{C \subset V \, | \, C$ is an irreducible curve
with $g^s(C) = C$ for some $s > 0 \}$ the set of all $g$-periodic curves.
For a divisor $M$ and subgroup $G \le \Aut(V)$, \,
$\Stab_G(M)$ $= \{x \in G \, | \, x^*M \equiv \alpha_x M$ in $\NS_{\C}(V)$
for some $\alpha_x \in \C\}$ is the {\it 'stabilizer' subgroup}.
For a group $G$, \, $Z(G) := \{g \in G \, | \, xg = gx$ for all $x \in G\}$
is the {\it centre} of $G$. If $H$ is a subgroup (resp. normal subgroup)
of $G$, we denote $H \le G$ (resp. $H \trianglelefteq G$).

\par
A {\it $(-n)$-curve} $C$ on a surface is a
curve with $C \cong \P^1$ and $C^2 = -n$.
A connected divisor on a surface is a {\it rational tree}, if it
is of simple normal crossing,
has the dual graph a tree
and consists of smooth rational curves.
\par \vskip 0.5pc
\begin{definition}\label{min}
Let $X$ be a smooth projective surface and $G \le \Aut(X)$ a
subgroup.
The pair $(X, G)$ is {\bf not} {\it minimal} (resp. is {\it minimal})
if the equivalent conditions below
are satisfied (resp. if neither of the conditions below is satisfied):
\begin{itemize}
\item[(1)]
There is a non-empty finite set $\Sigma$ of disjoint $(-1)$-curves on $X$
such that $\Sigma$ is $G$-stable, i.e., $G$ acts on $\Sigma$ as
permutations.
\item[(2)]
There is a $G$-equivariant non-isomorphic smooth blowdown $X \rightarrow Y$
onto a surface $Y$ endowed with a faithful $G$ action.
\end{itemize}
\end{definition}
\par \vskip 0.5pc
We use the format in KMM \cite{KMM} Theorem 7-3-1
for Fujita's result in \cite{Fu}:
\par \vskip 0.5pc
\begin{lemma}\label{Z-decomp}
Let $X$ be a smooth projective surface and $D$
a pseudo effective $\R$-divisor. Then there is a unique effective
$\R$-divisor $N$ satisfying:
\begin{itemize}
\item[(1)]
Either $N = 0$, or the irreducible components $N_i$
of $N$ give rise to a negative definite intersection matrix
$(N_i . N_j)_{1 \le i, j \le s}$;
\item[(2)]
$P : = D - N$ is nef; and
\item[(3)]
$P . N = 0$ (equivalently (using $(2)$), $P . N_i = 0$ for all $i$).
\end{itemize}
Finally $N$ ard $P$ are $\Q$-Cartier divisors if so is $D$.
\end{lemma}
\par \vskip 0.5pc
For $D$ in the lemma above, we write $D = P + N$ and call it
the {\it Zariski decomposition} for pseudo effective divisor $D$.
The uniqueness above also shows: if $D \equiv D'$
with $D = P + N$ and $D' = P' + N'$ their Zariski decompositions,
then $N = N'$ (equal) and $P \equiv P'$.
\par
For the result below, we may refer to BHPV \cite{BHPV}, Ch IV, Cor. 7.2
and the Cauchy-Schwartz inequality; see Kawaguchi \cite{Kawaguchi}
Lemma 1.2.
\par \vskip 0.5pc
\begin{lemma}\label{nef}
Let $X$ be a smooth projective surface. Let $M$ and $D$ be
$\R$-divisors on $X$ neither of which is numerically trivial.
Suppose that $M$ is nef, $D$ is pseudo-effective
and $M . D = 0$. The we have:
\begin{itemize}
\item[(1)]
$D^2 \le 0$; and
\item[(2)]
$D^2 = 0$ holds if and only if $M^2 = 0$ and $D \equiv a M$ for some $a > 0$.
\end{itemize}
\end{lemma}
\par \vskip 0.5pc
\begin{definition}\label{entropy} {\bf Topological entropy.}
Let $M$ be a compact K\"ahler manifold.
Let $\lambda(g)$ be the {\it spectral radius}
of the action $g^*$ on the cohomology ring $H^*(M, \C)$, i.e., the
maximum of moduli of its eigenvalues.
The {\it topological entropy} of $g$ is defined as
$e(g) = \log \lambda(g)$. See \cite{Gr}, \cite{Yo}
\cite{Fr}, also \cite{DS2}.
\par
It is known that $e(g) \ge 0$, and $e(g) = 0$ if and only if every
eigenvalue of the action $g^*$ above has modulus 1.
Also $e(g) > 0$ if and only if at least one
eigenvalue of the restriction $g^*|H^{1,1}(M)$ has modulus different from 1;
see Dinh and Sibony \cite{DS1}.
\par
An element $g$ is of {\it null} (resp. {\it positive}) entropy if $e(g) = 0$
(resp. $e(g) > 0$).
A subgroup $G \le \Aut(M)$ is of {\it null entropy} (resp. {\it positive entropy})
if $e(g) = 0$ for all $g \in G$ (resp. $e(g) > 0$ for at least one $g \in G$).
\par
When $\bar{M}$ is a normal projective surface and $g \in G \le
\Aut(\bar{M})$, we say $g$ (resp. $G$) is of null or positive entropy
if so is $g$ (resp. $G$) regarded as an element (resp. a subgroup) of $\Aut(M)$.
Here $M \rightarrow \bar{M}$ is Hironaka's equivariant resolution
and induces the natural inclusion $\Aut(\bar{M}) \subseteq \Aut(M)$.
\end{definition}
\par \vskip 0.5pc
\begin{setup}\label{Assump} {\bf Assumption.}
\rm{
From now on till Proposition \ref{cyclo} (but except Lemma \ref{L1}),
we assume that $X$ is a smooth projective
surface with an (infinite) automorphism $g$ such that
$g^* | H^{1,1}(X)$ has an eigenvalue $\lambda$
with $|\lambda| > 1$. Namely, assume that $g$ is of positive entropy.
In Lemma \ref{lambda} below,
it turns out that such $\lambda$ with
modulus $>1$ is unique, and we can write $\lambda = \lambda(g)$.
}
\end{setup}
\par
When $X$ is only a normal projective surface with
$\tilde{X} \rightarrow X$ the minimal resolution,
we let $\lambda(g) = \lambda(g_{\tilde X})$;
here $g \in \Aut(X)$ induces $g_{\tilde X} \in \Aut(\tilde{X})$.
\par
We remark that this $\lambda$ is an algebraic number and is
either a {\it Pistol number} or a {\it Salem number}
according to the degree of $\lambda$ over $\Q$ (degree 2 or bigger).
See Salem \cite{Sa} and McMullen \cite{McP2} \S 2.
\par
For the proof of the result below and the current formulation of it, see
Cantat \cite{Ca2} Theorem 2.1.5, Dinh and Sibony \cite{DS1}
Theorem 2.1 and McMullen \cite{McK3} Theorem 3.2 and Corolalry 3.3,
or Kawaguchi \cite{Kawaguchi} Theorem 2.1.
\par \vskip 0.5pc
\begin{lemma} \label{lambda} 
Let $X, g$ be as in $\ref{Assump}$. Then $\lambda = \lambda(g) > 1$;
$\lambda$ and $\lambda^{-1}$ are conjugate over $\Q$. Further,
setting $h = h^{1,1}(X)$, the following are all eigenvalues of $g^*
| H^{1,1}(X)$  (with each $|\alpha_j| = 1$; in particular,
$\lambda(g^{-1}) = \lambda(g)$):
$$\lambda, \, \lambda^{-1}, \, \alpha_1, \, \alpha_2, \, \dots, \, \alpha_{h - 2}.
$$
\end{lemma}
\par \vskip 0.5pc
For the proof of the result below, see Cantat \cite{Ca1} Theorem 2,
Diller and Favre \cite{DF} Theorem 5.1, or Kawaguchi \cite{Kawaguchi}
Proposition 2.5 and Lemma 3.8.
\par \vskip 0.5pc
\begin{lemma}\label{L1}
Let $X$ be a normal projective surface and $g$ an automorphism
with positive entropy $\log(\lambda)$.
Then we have:
\begin{itemize}
\item[(1)]
There exist nef $\R$-Cartier divisors $L^+ = L(g)^{+}$ and $L^- = L(g)^-$
(unique upto to positive scalars),
which are not numerically trivial, such that the following are true
(especially, $L(g^{\pm})^+ = L(g^{\mp})^-$):
$$g^*L^+ \equiv \lambda L^+, \hskip 1pc g^*L^- \equiv \lambda^{-1} L^-.$$
\item[(2)]
Let $\sigma : \tilde{X} \rightarrow X$ be the minimal
resolution so that $g \in \Aut(X)$ induces $g_{\tilde X} \in \Aut(\tilde{X})$.
Then $L(g_{\tilde X})^{\pm}$ equals $\sigma^*L(g)^{\pm}$ up to positive scalars.
\end{itemize}
\end{lemma}
\par \vskip 0.5pc
Indeed, by the proof
in Kawaguchi, Lemma \ref{nef} and reducing to $\tilde{X}$
(see also Cantat \cite{Ca1} Theorem 2 and Diller and Favre \cite{DF} Theorem 5.1),
we can start with any ample (or even big $\R$-) Cartier divisor $B$ on $X$ and take
$L^+$ and $L^-$ as follows:
$$L^+ = \lim_{n \rightarrow +\infty} \frac{(g^*)^n(B)}{\lambda^n}, \hskip 0.5pc
L^- = \lim_{n \rightarrow +\infty} \frac{(g^*)^{-n}(B)}{\lambda^n}.$$
\par
By the description above, we have the following
(see Kawaguchi \cite{Kawaguchi} Proposition 2.5; use also the Hodge index theorem
for (3) - (4)):
\par \vskip 0.5pc
\begin{lemma}\label{L2}
Let $X, g$ be as in $\ref{Assump}$. The following are true.
\begin{itemize}
\item
[(1)] $(L^+)^2 = 0 = (L^-)^2$.
\item
[(2)] $L: = L^+ + L^-$ is nef and big (i.e., $L^2 > 0$).
\item
[(3)] Suppose that $D$ is an $\R$-divisor on $X$ such that
$(g^s)^*D \equiv D$ for some $s > 0$. Then $L^+ . D = 0 = L^- . D$,
so either $D^2 < 0$ or $D \equiv 0$.
\item
[(4)] $L^+ . K_X = 0 = L^- . K_X$; so either $K_X^2 < 0$ or $K_X \equiv 0$.
\end{itemize}
\end{lemma}
\par \vskip 0.5pc
Here is the relation between $\Stab(g)$ and $\Null(L^{\pm})$:
\par \vskip 0.5pc
\begin{lemma}\label{L3}
Let $X, g$ be as in $\ref{Assump}$. The following are true.
\begin{itemize}
\item[(1)]
$\Stab(g) = \Null(L^+) \cap \Null(L^-) = \Null(L)$,
where $L = L^{+} + L^-$.
\item
[(2)]
$\Null(L)$ is either empty, or
a finite set with negative intersection matrix.
In particular, $|\Null(L)| < \rho(X)$, the Picard number.
\item
[(3)] If $\Null(M)$ is a finite set for $M = L^+$ or
$M = L^-$ (this is always true by Theorem $\ref{infmany}$),
then $\Null(M) = \Stab(g)$.
\item
[(4)] The pair $(X, \langle g \rangle)$ is minimal
if and only if $\Stab(g)$ does not contain any $(-1)$-curve.
\end{itemize}
\end{lemma}
\begin{proof}
(2) is the consequence of the Hodge index theorem and that $L^2 > 0$.
The nefness of $L^{\pm}$ implies the second
equality in (1). Clearly, $g$ stabilizes each of the sets
$\Null(L^+)$, $\Null(L^-)$ and $\Null(L)$.
So (1) and (3) follow (see (2) and Lemma \ref{L2};
see Kawaguchi \cite{Kawaguchi} Proposition 3.1).
The set $\Sigma$ (possibly empty) consisting of all $(-1)$-curves in $\Stab(g)$,
is $g$-stable. So (4) is just the definition in
\ref{min}.
\end{proof}
\par
We prove several frequently-used properties of $g$ of positive entropy.
The assertion (4) below follows from the proof of Kawaguchi \cite{Kawaguchi} Claim 3.8.1.
Denote by $P_g(x)$ the {\it characteristic} polynomial
of $g^* | \NS_{\Q}(X)$, and $f_g(x)$ the {\it minimal} (irreducible) polynomial
of $\lambda(g)$ over $\Q$.
\par \vskip 0.5pc
\begin{proposition}\label{cyclo} Let $X, g, \lambda, L^{\pm}$
be as in $\ref{Assump}$. The following are true.
\begin{itemize}
\item[(1)]
For $r \ge 1$ one has
$P_{g^r}(x) = f_{g^r}(x) q_{g^r}(x)$ where $q_{g^r}(x)$ is a product of cyclotomic polynomials while $f_{g^r}(x)$ has no root of $1$.
One has $q_{g^s}(x) = (x-1)^e$
for some positive integers $s = s(g)$ and $e = e(g)$.
Also $\deg(q_{g^r}(x)) = e$ for all $r \ge 1$.
\item[(2)]
Further, we have the decomposition into two $g$-stable $\Q$-spaces
$$\NS_{\Q}(X) = V_{g} \oplus V_{g^s = 1}$$
such that $(g^r)^* | V_g$ ($r \ge 1$) has $f_{g^r}(x)$ as its characteristic polynomial
and also minimal polynomial,
the inclusion (*) : $V_g \otimes_{\Q} \R \supset \R [L^+] \oplus \R [L^-]$,
$e = \rank_{\Q} V_{g^s = 1}$, and
$$V_{g^s = 1} = \{v \in \NS_{\Q}(X) \, | \, g^s(v) = v\}
= (L^+)^{\perp} \cap (L^-)^{\perp} \cap \NS_{\Q}(X).$$
\item[(3)]
If $X$ is rational then
$\Ker(\Aut(X) \rightarrow \Aut(\Pic(X)))$ is finite.
\item[(4)]
We may take $L^{\pm}$ to be in $(\Pic(X)) \otimes_{\Z} \Z[\lambda]$.
\item[(5)]
$\lambda$ is not a rational number.
\item[(6)] None of positive multiples of
$L^{\pm}$ is $\Q$-Cartier.
\end{itemize}
\end{proposition}
\begin{proof}
(1) The first part follows from Lemma \ref{lambda} and Kronecker's theorem;
see McMullen \cite{McK3} Corollary 3.3.
For the second, let $s$ be the lcm of orders of the
roots in $q_g(x)$.
\par
(2) By (1), we have $\NS_{\Q}(X) = V_r' \oplus V_r''$
(each summand being $g$-stable)
so that $f_{g^r}(x)$ and $q_{g^r}(x)$ are respectively
the characteristic polynomials of
$(g^r)^* | V_r'$ and $(g^r)^* | V_r''$. Further, we have $V_i' = V_j'$
(denoted as $V_g$)
and $V_i'' = V_j''$ (denoted as $V_{g^s=1}$) for all $i, j \ge 1$,
since the only $g$-stable vector subspaces of $V_r'$
are $\{0\}$ and itself by the irreducibility of $f_{g^r}(x)$.
\par
Since the classes $[L^{\pm}]$
are eigenvectors w.r.t. eigenvalues $\lambda^{\pm}$ of $g^*$,
we have the inclusion (*) in (2).
\par
For the second equality in the second display of (2),
in view of Lemma \ref{L2},
we have only to show the assertion that
$(g^r)^* | W  = \id$ for some $r \ge 1$,
where $W: = (L^+)^{\perp} \cap (L^-)^{\perp} \cap \NS(X)$.
By the Hodge index theorem, $W$ (modulo its finite torsion)
is a negative definite integral lattice.
So $\Aut(W)$ is finite and the assertion is true. Thus the second equality is proved.
\par
Now the first equality (even when $s$ is replaced by $s n$ with $n \ge 1$)
follows from the second, the choice of $s$ and the application
of Lemma \ref{L2} (3) inductively on the size of
the Jordan canonical form of $g^* | V_r'' \otimes_{\Q} \C$.
\par
(3) Suppose the contrary that this Kernel is infinite. Then $X$ has finitely many
'exceptional curves' by Harbourne \cite{HPAMS} Proposition 1.3 and its terminology
at the last paragraph of page 409. Since $g$ acts on the finite set of these
curves, for some common $m \ge 1$,
each irreducible component of these exceptional curves
is stabilized by $g^m$.
Let $X \rightarrow Y$ be the $\langle g^m \rangle$-equivariant
smooth blowdown to a relatively minimal model. Then $K_Y^2 \ge 8$,
while $K_Y^2 < 0$ by Lemma \ref{L2} for $g^m$ being of positive entropy on $Y$.
It is absurd.
\par
(5) If $\lambda \in \Q$, then $\lambda \in \Z$
because it is algebraic over $\Q$, so $x - \lambda$ is its
minimal polynomial over $\Q$, contradicting the fact that
$\lambda$ and $\lambda^{-1}$ are conjugate over $\Q$ in
Lemma \ref{lambda}.
\par
(6) If $L^+$ is $\Q$-Cartier say, then
by intersecting $g^*L^+ \equiv \lambda L^+$ with a
Cartier ample divisor $H$, we see that
$\lambda = (H . g^*L^+)/(H . L^+)$ is a rational number.
This contradicts (5).
The proposition is proved.
\end{proof}
\par
Here is another consequence of the uniqueness result of Lemma \ref{L1}.
\begin{lemma}\label{Lcor}
Let $X$ be a smooth projective surface, $M$ a nef $\R$-divisor which
is not numerically trivial, and $g \in \Aut(X)$ such that
$g^*M \equiv \alpha M$ for some $\alpha \in \C$. Then we have:
\begin{itemize}
\item[(1)]
Either
$\alpha = 1$, or $\alpha > 1$, or $0 < \alpha < 1$.
\item[(2)]
If $\alpha = 1$, then
$g$ is of null entropy.
\item[(3)]
If $\alpha^{\pm} > 1$, then
$g$ is of positive entropy with $\lambda(g) = \alpha^{\pm}$ and
$M$ is equal to $L(g)^{\pm}$ (up to a positive scalar).
\end{itemize}
\end{lemma}
\begin{proof}
Intersecting the equality $g^*M \equiv \alpha M$ with an ample divisor $H$, we
get $\alpha = (H . g^*M)/(H . M) > 0$. Then the lemma follows from
Lemmas \ref{lambda} - \ref{L2} and that $M^2 \ge 0$ for $M$ being nef.
\end{proof}
The result below shows that when dealing with automorphisms of
positive entropy we may quotient away a finite group.
\par
Let $\sigma : \Aut(X) \rightarrow \Aut(\NS(X))$
be the natural homomorphism. For $g \in G \le \Aut(X)$,
set $\tilde{g} = \sigma(g)$ and $\tilde{G} = \sigma(G)$.
\par \vskip 0.5pc
\begin{lemma}\label{quot}
Let $X$ be a smooth projective surface and $H \le \Aut(X)$
a subgroup such that $\sigma(H) = \tilde{H}$ is a finite subgroup of $\Aut(\NS(X))$.
Suppose that $\sigma(g) = \tilde{g} \in N_{\Aut(\NS(X))}(\tilde{H})$
(i.e., $\tilde{g} \tilde{H} \tilde{g}^{-1} = \tilde{H}$).
Then we have:
\begin{itemize}
\item[(1)]
If $g$ is of positive entropy, then for both $M = L(g)^+$ and $M = L(g)^-$ we have
$h^*M \equiv M$ for all $h \in H$.
\item[(2)]
Suppose that $H$ is already finite.
Let $\pi : X \rightarrow \bar{X} = X/H$ be the
quotient map. Then $g$ is of positive entropy if and only if so is
$\bar{g} = gH \in \bar{G} = G/H \le \Aut(\bar{X})$.
If this is the case, we have $\lambda(g) = \lambda(\bar{g})$
and can take $L(g)^{\pm} = \pi^*L(\bar{g})^{\pm}$.
\end{itemize}
\end{lemma}
\begin{proof}
(1) Consider the case $M = L(g)^+$ (the other case is similar),
so $g^*M \equiv \lambda M$ with $\lambda = \lambda(g) > 1$.
Set $M' = \sum_{x \in \tilde{H}} x^*(M)$ (identifying
$M$ with its class in $\NS(X) \otimes_{\Z} \R$).
Then we have (where $h \in H$):
$$g^*(M') = \sum_{x \in \tilde{H}}
(\tilde{g}^{-1} x \tilde{g})^* \tilde{g}^* M \equiv
\lambda \sum_{y \in \tilde{H}} y^*M
= \lambda M',$$
$$h^*(M') = \sum_{x \in \tilde{H}} (x \tilde{h})^*M = \sum_{y \in \tilde{H}} y^*M
= M'.$$
Clearly, $M'$ is nef and is not numerically trivial,
so the uniqueness Lemma \ref{L1} implies that $M \equiv \alpha M'$ for some $\alpha > 0$.
Thus (1) follows.
\par
(2) Suppose that $\bar{g}$ is of positive entropy.
So $\bar{g}^*L(\bar{g})^+ \equiv \lambda L(\bar{g})^+$
with $\lambda = \lambda(\bar{g}) > 1$.
Set $L(g)^+ = \pi^*L(\bar{g})^+$ which is nef and is not
numerically trivial. Then $g^*L(g)^+ \equiv \lambda L(g)^+$ because
$\pi \circ g = \bar{g} \circ \pi$. Thus $g$ is of positive
entropy $\log \lambda(g) = \log \lambda$.
\par
Conversely, suppose that $g$ is of positive entropy $\log \ \lambda = \log \lambda(g)$.
We shall use the fact below about the $H$-invariant sublattice
$$((\Pic(X))^H \otimes_{\Z} \R = \pi^*(\Pic(\bar{X})) \otimes_{\Z} \R.$$
By the proof of (1), $L(g)^+$ (re-chosen like $M'$ above)
belongs to the LHS above, so it is also in the RHS.
Thus $L(g)^+ = \pi^* \bar{L}$ for some $\R$-Cartier divisor $\bar{L}$.
Since $\pi \circ g = \bar{g} \circ \pi$, we have
$\pi^*(\lambda \bar{L}) \equiv \pi^*(\bar{g}^*\bar{L})$.
Hence $\bar{g}^*\bar{L} \equiv \lambda \bar{L}$ by the injectivity
of $\pi^*: (\Pic(\bar{X})) \otimes_{\Z} \R \rightarrow
(\Pic(X)) \otimes_{\Z} \R$. So $\bar{g}$ is of positive entropy.
Thus (2) and the lemma are proved, since the case $L(g)^-$ is similar.
\end{proof}
\section{Groups stabilizing nef classes; proof of Theorem \ref{Thsv}}
In this section, we will determine the dynamical structure of
the stabilizer subgroup $\Stab_G(M)$
when $M$ is nef and is not numerically trivial.
\par
We will also determine the dynamical group structure of the full
automorphism group $\Aut(S^v)$ for McMullen's surfaces $S^v$ in \cite{McP2} \S 7.
\par
The following result is very important in proving our main theorems.
\par \vskip 0.5pc
\begin{theorem}\label{nefstab1}
Let $X$ be a smooth projective surface and $G \le \Aut(X)$ a
subgroup. Let $M$ be a nef $\R$-divisor which is not numerically
trivial. Suppose that $G = \Stab_G(M)$ and $G$ is of positive
entropy. Then we have:
\begin{itemize}
\item[(1)] One has $G = \langle h_m \rangle \ltimes T$ (semi-direct product)
with $h_m$ of positive entropy (and $M = L(h_m)^+$) and
$\Imm (T \rightarrow \Aut(\NS(X)))$ finite.
\item[(2)] For $g \in G$, one has $g \in T$ if and only if $g$ is of null entropy;
if and only if $g^*M \equiv M$ holds.
\item[(3)] Set $h:= h_m$. Then $T | V_{h} = \id$. Here $\NS_{\Q}(X)
= V_{h} \oplus V_{h^{s(h)}=1}$ is as in $\ref{cyclo}$.
\end{itemize}
\end{theorem}
\begin{proof}
By the assumption, $G$ contains an automorphism $g$ of positive
entropy $\log \lambda(g) > 0$. By Lemma \ref{Lcor} and switching $g$
with $g^{-1}$ if necessary, we may assume that $M = L(g)^+$ and
$g^*M \equiv \lambda(g) M$. Also, in the definition of $\Stab_G(M)$,
either $\alpha_x = 1$ and $x$ is of null entropy, or $\alpha_x^{\pm}
> 1$ and $x$ is of positive entropy $\lambda(x) = \alpha_x^{\pm}$.
\par
Consider the map $\varphi: G \rightarrow \R$ which takes $x$ to $\log\alpha_x$
if $x^* M \equiv \alpha_x M$.
Clearly, $\varphi$ is a homomorphism onto an abelian subgroup of the (torsion-free) additive group $\R$.
\par
We assert that $\Imm(\varphi)$ does not contain a subgroup of the type $\Z \times \Z$.
Suppose the contrary that $\langle \log(\alpha_{g_1}) \rangle \times \langle \log(\alpha_{g_2}) \rangle
\cong \Z \times \Z$ is a subgroup of $\Imm(\varphi)$. Here
$g_i \in G$ with $g_i^*M \equiv \alpha_{g_i} M$.
Thus $1, \log(\alpha_{g_2})/\log(\alpha_{g_1})$ are linearly independent over $\Q$.
Now by the classical Kronecker's Theorem (or Dirichlet's Theorem),
for any $\varepsilon > 0$, there are integers $m_i$ such that:
$$|m_2 \frac{\log(\alpha_{g_2})}{\log(\alpha_{g_1})} - m_1| < \varepsilon.$$
After relabelling, we may assume that:
$$0 < m_1 \log(\alpha_{g_1}) - m_2 \log(\alpha_{g_2}) < \varepsilon_1 := \varepsilon \max\{|\log(\alpha_{g_1})|, |\log(\alpha_{g_2})|\}.$$
Thus for $h = g_1^{m_1}/g_2^{m_2} \in G$,
we have $h^*M \equiv \alpha_h M$ where
$\alpha_h := \alpha_{g_1}^{m_1}/\alpha_{g_2}^{m_2}$
satisfies $1 < \alpha_h < e^{\varepsilon_1}$. So $h$ is of positive entropy
with $\lambda(h) = \alpha_h$.
On the other hand,
in \cite{McP2} Theorem 1.2, McMullen has proved that $\lambda(f) \ge \lambda(\text{\rm Lehmer})
\thickapprox 1.17628081$ for every surface automorphism $f$ of positive entropy.
We get a contradiction if we let $\varepsilon$ (and hence $\varepsilon_1$) tend to zero.
This proves the assertion.
\par
Since $\varphi$ factors through $G \subseteq \Aut(X) \rightarrow \Aut(\NS(X))$
while the latter is countable, we can write $\Imm(\varphi) = \cup_{n \ge 1} \bar{G}_n$
where each $\bar{G}_n \le \R$ is finitely generated (and abelian).
By the assertion above and noting that $\varphi(g) = \log \lambda(g) > 0$
for some $g \in G$
(and the fundamental theorem for f.g. abelian groups)
we may assume that $\bar{G}_n = \langle \log \lambda(g_n) \rangle$
for some $g_n \in G$ of positive entropy $\log \lambda(g_n)$.
Now $\bar{G}_n \le \bar{G}_{n+1}$ implies, by induction, that
$\lambda(g_n) = (\lambda(g_1))^{1/s_n}$ for some positive integer $s_n$.
Applying McMullen's result above again, we have
$(\lambda(g_1))^{1/s_n} = \lambda(g_n) \ge \lambda(\text{\rm Lehmer}) > 1.1$
for all $n \ge 1$.
Thus there is a constant $N$ such that
$s_N = s_{N+1} = \dots$,
whence $\Imm (\varphi) = \bar{G}_N = \langle \log \lambda(h_m) \rangle$,
where we set $h_m := g_N$ and may assume that $M = L(h_m)^+$.
\par
Set $T := \Ker(\varphi)$. Then $T$ is of null entropy;
so the assertion (2) follows (see also Lemma \ref{Lcor}).
Since $h_m$ is of positive entropy, $\langle h_m \rangle \cap T = (1)$
and hence $G = \langle h_m \rangle \ltimes T$.
This proves (1) (except the finiteness of $T | \NS(X)$).
\par
(3) Write $h := h_m$ for simplicity.
We use the notation in Proposition \ref{cyclo}.
Note that $f(x) := f_h(x)$ is the minimal polynomial of $\lambda = \lambda(h)$
over $\Q$ (hence $f(x)$ has only simple zeros by the Galois theory) and also the
minimal polynomial of $h^* | V_h$,
whence $V_h \otimes_{\Z} K(f)$ is spanned by the eigenvectors
$v_{\beta}$ w.r.t.
eigenvalues $\beta$ for $h^* | V_h \otimes_{\Z} K(f)$.
Here $K(f)$ is the splitting field ($\subset \C$)
of $f(x)$ over $\Q$.
\par
So we have only to show the assertion that
$t^*v_{\beta} \equiv v_{\beta}$ for all such
eigenvector $v_{\beta}$
(so that $nt^*v_{\beta}$ and $nv_{\beta}$ are algebraically equivalent for some $n > 0$
and hence $t^*v_{\beta}$ equals $v_{\beta}$ in $\NS_{\Q}(X)$).
There is a $\gamma = \gamma_h$ in the Galois group $\Gal(f) = \Gal(K(f)/\Q)$ of $f(x)$
such that $\gamma^{-1}(\beta)$ equals the spectral radius $\lambda$
(see Lemma \ref{lambda}).
Set $L^{\pm} = L(h)^{\pm} \in (\Pic(X)) \otimes_{\Z} \Z[\lambda]$
(see Proposition \ref{cyclo}),
so that $h^*L^{\pm} \equiv \lambda^{\pm} L^{\pm}$.
For $g \in \Aut(X)$, extend $g^* | \Pic(X)$ \, $\C$-linearly to
$g^* | (\Pic(X)) \otimes_{\Z} \C$.
For a $\C$-divisor $D = \sum c_i D_i$ with $D_i \in \Pic(X)$
and $c_i \in K(f)$,
define $\gamma^*(D) = \sum \gamma(c_i) D_i$.
Then $\gamma \circ g = g \circ \gamma$.
Now:
$$h^*\gamma^*(L^+) = \gamma^* h^*(L^+) \equiv \gamma^*(\lambda L^+) = \beta \gamma^*(L^+).$$
So $v_{\beta} := \gamma^*(L^+)$ (its class, to be precise) is an eigenvector w.r.t. $\beta$
for $h^* | \NS(X) \otimes_{\Z} K(f)$.
For $t \in T$, we have $t^*v_{\beta} = t^*\gamma^*(L^+) = \gamma^* t^*(L^+)
\equiv \gamma^*(L^+) = v_{\beta}$. The assertion and hence (3) are proved.
\par
Finally, to show $T | \NS(X)$ is finite it is enough to show that $T | \NS_{\Q}(X)$ is finite, since $\NS(X)$ has a finite torsion.
By (3) and the proof of Proposition \ref{cyclo}, we then have only to show
$T$ stabilizes $W: = (L^+)^{\perp} \cap (L^-)^{\perp} \cap \NS(X)$
(noting that $W \otimes_{\Z} \Q = V_{h^{s(h)}=1}$ and $O(W)$ is finite).
This is clear because $t^*L^{\pm} \equiv L^{\pm}$ for $t \in T$
by (3) and Proposition \ref{cyclo} (2). We are done.
\end{proof}
By the proof above, $\langle h_m \rangle \ltimes T$ stabilizes (up to scalars)
$L(h_m)^-$. Reversing the process or noting that $L(h_m)^- = L(h_m^{-1})^+$, one obtains:
\par \vskip 0.5pc
\begin{corollary}
Let $X$ be a smooth projective surface and $g$ an automorphism
of positive entropy. Then $\Stab_{\Aut(X)}(L(g)^+) = \Stab_{\Aut(X)}(L(g)^-)$.
\end{corollary}
We now elaborate further about
the stabilizer subgroup $\Stab_G(M)$.
\par \vskip 0.5pc
\begin{proposition}\label{nefstab2}
Assume that $G = \Stab_G(M)$ and $G$ is of positive entropy as in
$\ref{nefstab1}$. Assume further that $G \trianglelefteq F \le
\Aut(X)$. Then we have:
\begin{itemize}
\item[(1)]
One has $|F : \Stab_F(M)| \le 2$, so
$\Stab_F(M)$ is normal in $F$. Write $\Stab_F(M) = \langle h_F \rangle
\ltimes T_F$ with $M = L(h_F)^+$, as in $\ref{nefstab1}$.
\item[(2)]
Suppose that $|F : \Stab_F(M)| = 2$ and let $\tau \in F \setminus \Stab_F(M)$.
Then $F = \{\tau^i h_F^j t \, | \, i = 0,1; \ j \in \Z; \ t \in T_F\}$.
Also $g \in F$ is of positive entropy if and only if $g = h_F^j t$
for some $j \ne 0$ and $t \in T_F$.
\end{itemize}
\end{proposition}
\begin{proof}
(1) We use Theorem \ref{nefstab1} (and its notation) :
$G = \langle h_m \rangle \ltimes T$.
Write $h = h_m$ for simplicity. Take $f \in F$.
Then $f h f^{-1} = h^r t$ for some $r = r(h) \in \Z$ and $t \in T$.
Applying the equality to $M = L(h)^+$, we get $(f h f^{-1})^* M \equiv \lambda^r M$
(with $\lambda = \lambda(h) > 1$)
and $h^*(f^*M) \equiv \lambda^r (f^*M)$. By the uniqueness Lemma \ref{L1},
either $r = 1$ and $f^*M$ equals $L(h)^+$, up to a positive scalar
(so $f \in \Stab_F(M)$),
or $r = -1$ and $f^*M$ equals $L(h)^-$,
up to a positive scalar.
\par
Suppose there are $f_i \in F$ such that $f_i h f_i^{-1} = h^{r_i} t_i$
as above, but with $r_1 = r_2 = -1$.
Then $(f_2 f_1) h (f_2 f_1)^{-1} = f_2(h^{-1} t_1) f_2^{-1}
= (h^{-1} t_1)^{-1} t_1' = h t'' t' =: h t'''$.
Here $t_1' = f_2 t_1 f_2^{-1}$ and $t'''$ are all in $T$.
Hence, by the argument above, $f_1 f_2$ (and especially $f_i^2$)
are all in $\Stab_F(M)$. Thus (1) follows. Indeed,
we may assume that $h_m = h_F^{\ell}$ for some $\ell \ge 1$,
and $L(h_F)^{\pm} = L(h)^{\pm}$.
\par
(2) The first part is from (1). Set $L^{\pm} := L(h_F)^{\pm}$.
We have $\tau T_F \tau^{-1} = T_F$
since $T_F$ is the set of all null entropy elements in $\Stab_F(M)$.
Since $h_F$ is of positive entropy,
by the proof of (1), we have $\tau h_F \tau^{-1} = h_F^{-1} t$ for some $t \in T_F$
and may assume that $\tau^*L^+ = L^-$ and ($\tau^*L^- =$) \,
$(\tau^2)^*L^+ \equiv \alpha L^+$
for some $\alpha > 0$.
We also have $(\tau^2)^* L^- \equiv \tau^*(\alpha L^+)
= \alpha L^-$.
By Lemma \ref{Lcor} and the uniqueness Lemma \ref{L1},
we have $\alpha = 1$, and hence $\tau^2$ (and also $\tau$) are of null entropy.
Thus $\tau^2 = t' \in T_F$. Now for $t_F \in T_F$, we see that
$(\tau h_F^j t_F)^2 = (\tau h_F^j t_F \tau^{-1}) t' h_F^j t_F
= (h_F^{-1} t)^j t'' h_F^j t_F \in T_F$ (here $t'' = (\tau t_F \tau^{-1}) t' \in T_F$),
whence $\tau h_F^j t_F$ is of null entropy. (2) follows.
The proposition is proved.
\end{proof}
Let
$\sigma : \Aut(X) \rightarrow \Aut(\NS(X))$
be the natural homomorphism. We have:
\par \vskip 0.5pc
\begin{proposition}\label{nefstab3}
Let $X$ be a smooth projective surface and
$G \le \Aut(X)$ a subgroup of positive entropy.
Suppose that there is an $H \triangleleft G$ such
that $G/H$ is abelian and $\sigma(H)$
is a finite subgroup of $\Aut(\NS(X))$.
Then $G = \Stab_G(P)$
for every $P = L(g)^{\pm}$ with $g \in G$ being positive entropy.
\end{proposition}
\begin{proof}
By Lemma \ref{quot}, we have $h^*P \equiv P$ for every $h \in H$.
By the assumption, for every $g_1 \in G$
we have $g_1 g = h g g_1$ for some $h \in H$.
Applying this equality to $P = L(g)^+$, we get
$g^*(g_1^*P) \equiv \lambda (g_1^*P)$ with $\lambda = \lambda(g) > 1$.
By the uniqueness Lemma \ref{L1},
$g_1^*P$ equals $L(g)^+$ (up to a positive scalar). So $g_1 \in \Stab_G(P)$.
Thus $G = \Stab_G(P)$. The case $P = L(g)^-$ is similar.
The proposition is proved.
\end{proof}
In \cite{McP2} \S 7, McMullen constructed rational surfaces $S^v$ for all
$v \in (\C^{1,n})^*$ there. The latter set contains all
leading eigenvectors $v$ w.r.t. an eigenvalue $\beta$,
for a Coxeter element $w$
in the Weyl group $W_n$ ($n \ge 10$)
of the Minkowski parabolic lattice of rank $n+1$
(so $S^v$ has Picard number $n+1$).
Here $v$ is a {\it leading eigenvector} for $w$ if
$\beta$ is conjugate over $\Q$ to the largest
(in terms of modulus) eigenvalue $\lambda(w)$ of $w$.
One has $\lambda(w) > 1$ (see \cite{McP2} before 2.5).
\par \vskip 0.5pc
\begin{setup} {\bf Proof of Theorem \ref{Thsv}}
\end{setup}
Set $S := S^v$.
Let $\sigma : \Aut(S) \rightarrow \Aut(\Pic(S))$ be the natural
homomorphism.
Since  $S$ is rational, we do the identification: $\NS(X) = \Pic(S) = \Z^{1,n}$ and
$(\Pic(S)) \otimes_{\Z} \C = \C^{1,n}$, as in \cite{McP2}.
Set $k_n = [K_S] \in \Pic(S)$. By \cite{McP2} Corollary 7.2, we have
$$\Aut(S) \cong W_n^v = \{w' \in W_n \, | \, [v]
\in \C^{1,n}/\C k_n \,\, \text{\rm is an eigenvector for} \,\, w' \}.$$
Since every $h \in \Aut(S)$ stabilizes $K_S$ and $K_S^{\perp}$, and $v \in K_S^{\perp}$,
we have $\sigma(\Aut(S)) \subseteq \Stab_{\Aut(\Pic(S))}(v)$.
So for every $h \in \Aut(S)$, we have $h^*(v) = \beta_h v$ for some $\beta_h \in \C$.
Fix any $g \in \Aut(S)$ (of positive entropy) realizing the Coxeter element $w$.
Set $\beta = \beta_g$.
Let $f(x) = f_g(x)$ be the minimal polynomial of $\beta_g$ over $\Q$.
Let $K(f_g)$, or $K(f_g, f_h)$ be the splitting field ($\subset \C$)
of $f_g(x)$, or $f_g(x)$ and $f_h(x)$
over $\Q$.
There is a $\tau = \tau_g$ in the Galois group $\Gal(f) = \Gal(K(f)/\Q)$ of $f(x)$
such that $\tau(\beta)$ equals the spectral radius $\lambda = \lambda(g) = \lambda(w)$
(see Lemma \ref{lambda}). Clearly, we may take $v \in (\Pic(S)) \otimes_{\Z} K(f)$.
Set $L^{\pm} = L(g)^{\pm} \in (\Pic(S)) \otimes_{\Z} \Z[\lambda]$
(see Proposition \ref{cyclo}),
so that $g^*L^{\pm} \equiv \lambda^{\pm} L^{\pm}$.
\par
We claim that $\Aut(S) = \Stab_{\Aut(S)}(M)$ for both $M = L^{\pm}$,
which will in turn imply the theorem by Theorem \ref{nefstab1}
and Proposition \ref{cyclo}.
\par
Take $h \in \Aut(S)$. Extend $h^* | \Pic(S)$ \, $\C$-linearly to
$h^* | (\Pic(S)) \otimes_{\Z} \C$.
Denote by the same $\tau$ its extension to an automorphism of
the splitting field $K(f_g, f_h)$ (by the isomorphism extension theorem).
For a $\C$-divisor $D = \sum c_i D_i$ with $D_i \in \Pic(S)$
and $c_i \in K(f_g, f_h)$,
define $\tau^*(D) = \sum \tau(c_i) D_i$.
Then $\tau \circ h = h \circ \tau$. We have:
$$\lambda \tau^*(v) = \tau^*(\beta v) = \tau^*(g^*(v))
= g^*(\tau^*(v)).$$
By Lemma \ref{lambda},
$\tau^*(v)$ equals the class of $\delta L^+$ for some $\delta \in K(f)$
(by the choice of $v$ and $L^+$).
Replacing $v$, we may assume that $\tau^*(v)$ equals the class of $L^+$.
Now the claim for $M = L^+$ (and hence the theorem) follows from the calculation below
(noting that the case $M = L^-$ is similar):
$$h^*L^+ = \tau^*h^*(\tau^{-1})^*L^+ \equiv \tau^*h^*(v) = \tau^* (\beta_h v)
= \tau(\beta_h) \tau^*(v) \equiv \tau(\beta_h) L^+.$$
This proves the theorem.
\section{Automorphisms of null entropy}
In this section, we consider groups of automorphisms of null entropy.
The following is the main result of the section.
\par
Let $\sigma : \Aut(X) \rightarrow \Aut(\Pic(X))$ be the natural homomorophism.
\par \vskip 0.5pc
\begin{theorem}\label{Thnulla}
Let $X$ be a smooth projective surface with irregularity
$q(X) = 0$ and let $G \le \Aut(X)$ be a subgroup of
null entropy such that $\sigma(G)$ is an infinte subgroup of
$\Aut(\Pic(X))$.
Suppose that the pair $(X, G)$ is minimal.
Then we have:
\begin{itemize}
\item [(1)]
There is a nef $\Q$-divisor $M$ (which might be zero) with $M^2 = 0$,
such that either $K_X \equiv M$ or $K_X \equiv -M$.
In particular, $K_X^2 = 0$.
\item[(2)]
If $X$ is rational, then $-K_X$ is nef and
$X$ has an anti $1$-canonical curve.
\end{itemize}
\end{theorem}
\par \vskip 0.5pc
\begin{setup}\label{setup4.1}
\rm{
We need some preparation before starting the proof.
For a smooth projective surface $X$,
let $N_1(X)$ be the $\R$-vector space (of rank $\rho(X)$, the Picard number)
of 1-cycles modulo numerical equivalence; let $\overline{\NE}(X)$
be the closure in $N_1(X)$ of the cone $\NE(X) \subset N_1(X)$
of effective 1-cycles (modulo numerical equivalence); see Kollar and Mori
\cite{KM}, Definitions 1.16 and 1.17. In our surface case,
an $\R$-divisor $D$ is pseudo-effective if and only if $D \in \overline{\NE}(X)$.
\par
Let $X$ be a smooth projective surface.
Set $L = \NS(X)/\tor$, the Neron-Severi lattice modulo its (finite) torsion.
Suppose that $G \le \Aut(X)$ is of null entropy.
Then the image $\overline{G} := \Imm(G \rightarrow \Aut(L))$
is of null entropy in the sense of Oguiso \cite{Og2}, \S 2.
Conversely, if $\overline{G}$ is of null entropy in the sense of
Oguiso, then $G$ is of null entropy in the usual sense of Definition \ref{entropy},
by virtue of Lemma \ref{L1}.
\par
We now assume that $q(X) = 0$. Then
$L = \Pic(X)/\tor$ and $\Pic(X) \cong L \oplus (\text{\rm finite torsion subgroup})$.
Clearly, $\sigma(G) \le \Aut(\Pic(X))$ is infinite if and only
if $\bar{G} \le \Aut(L)$ is infinite.
\par
For divisors $D_1, D_2$ {\it we say $D_1 = D_2$ in $L$}
if their classes in $L$ are identical, i.e., if $D_1 \equiv D_2$
(noting that $q(X) = 0$ is assumed).
}
\end{setup}
\par
The following result is proved in Oguiso \cite{Og2} Lemma 2.8.
We slightly change the formulation, where
$g(v)$ etc. should be understood as $\overline{g}(v)$
with $\overline{g}$ the image of $g$ via the surjective
homomorphism $G \rightarrow \overline{G}$.
\par \vskip 0.3pc
\begin{lemma}\label{Og2}
The following are true.
\begin{itemize}
\item[(1)]
There is a unique ray $0 \ne \R_{> 0} v \in \overline{\NE}(X)$ such that $v^2 = 0$
and $g(v) = v$ in $L$ for all $g \in G$. We can choose $v$ to be in $L$.
\item [(2)]
Suppose that $B$ is a subgroup of $O(L)$ and $\ell \in L$ such that
$\ell^2 > 0$ and $b(\ell) = \ell$ for all $b \in B$. Then $B$ is finite.
\end{itemize}
\end{lemma}
\par \vskip 0.5pc
We continue the proof of Theorem \ref{Thnulla}.
Clearly, (2) is a consequence of (1)
and the Riemann-Roch theorem. We now prove (1).
\par
Let $V$ be in $\Pic(X)$ whose class in $L$ is $v$.
We assert that $K_X . V = 0$. Suppose the contrary that
$K_X . V \ne 0$. Then there is a positive integer $t$ such that $D = K_X + tV$
satisfies $D^2 > 0$. Clearly, $g(D) = D$ in $L$ for all $g \in G$.
Then $\overline{G}$ is finite by Lemma \ref{Og2}, contradicting
the assumption on $\overline{G}$. So the assertion is true.
\par
Next we assert that $V$ is nef.
Let $V = P + N$ be the Zariski-decomposition.
Since $g(v) = v$ in $L$ (and hence $g(V) \equiv V$) for every $g \in G$
and by the uniqueness in Lemma \ref{Z-decomp},
we have $g(N) = N$ and $g(P) \equiv P$ (so $g(nP) = nP$ in $L$
for some positive integer $n$ with $nP$ integral)
for every $g \in G$. If $P^2 > 0$, then we get a contradiction as above.
So $P^2 = 0$. Now $0 = v^2 = (P+N)^2 = N^2$ implies that $N = 0$
and $V = P$ is nef. The assertion is proved.
\par
Serre duality says $h^2(X, K_X + V) = h^0(X, -V) = 0$ because $V . H > 0$
for an ample divisor $H$ on $X$. Since $q(X) = 0$,
one has $\chi(\OO_X) = 1 - q(X) + p_g(X) > 0$.
Now the Riemann-Roch theorem implies that $h^0(X, K_X + V) \ge \chi(\OO_X) > 0$.
If $K_X + V \equiv 0$, then $-nK_X \sim nV$ for some integer $n > 0$
because $q(X) = 0$. We are done.
\par
Suppose $K_X + V$ is not numerically trivial.
Let $K_X + V = P' + N'$ be the Zariski decomposition.
Since $g(K_X + V) \equiv K_X + V$ for all $g \in G$,
we have, as above, $g(P') = P'$ in $L$, \, $(P')^2 = 0$ and
$g(N') = N'$, whence $G$ permutes components of $N'$.
By the uniqueness of $v$, we have $P' \equiv a V$ for some $a \ge 0$.
So $K_X + (1-a) V \equiv N'$.
\par
If $N' = 0$, then $-n' K_X \sim n'(1-a) V$ for some integer $n' > 0$
and we are done.
Suppose that $N' \ne 0$.
Then $(N')^2 < 0$ and hence for some component $N_1$ of $N'$,
we have $0 > N_1 . N' = N_1 . (P' + N') = N_1 . (K_X + V) \ge N_1 . K_X$.
So $N_1$ is a $(-1)$-curve in $N'$. Let $\Sigma$ be the set of
$(-1)$-curves $N_j$ in $N'$. Then $N_i \cap N_j = \emptyset$ ($i \ne j$)
by the negativity of $N'$. Since $G(N') = N'$, we have
$G(\Sigma) = \Sigma$, contradicting the minimality of $(X, G)$.
The theorem is proved.
\section{Dynamics of groups; Proofs of some of Theorems \ref{Thnull} - \ref{Thtits}}
In this section we give dynamical structures for certain
groups of surface automorphisms.  We also prove Theorems \ref{Thnull}, \ref{Thz2} and \ref{Thtits}.
\par \vskip 0.5pc
\begin{setup} {\bf Proof of Theorem \ref{Thnull}.}
\end{setup}
\par
Clearly, there is a $G$-equivariant smooth
blowdown such that the pair $(Y, G)$ is minimal.
$G \le \Aut(Y)$ is also of null entropy.
So we may assume that $(X, G)$ is already minimal.
Indeed, $\Imm (G \rightarrow \Aut(\Pic(X)))$ is infinite
if and only if so is $\Imm (G \rightarrow \Aut(\Pic(Y)))$
since $G$ acts on the (finite) set of curves in the
exceptional divisor $E$ of $X \rightarrow Y$
and $\Pic(X)$ is the direct sum
of the pull back of $\Pic(Y)$ and the lattice generated by the
curves in $E$.
Set $\overline{G} = \Imm (G \rightarrow \Aut(\Pic(X)))$.
Consider the natural exact sequence:
$$1 \rightarrow J \rightarrow G \rightarrow \overline{G} \rightarrow 1.$$
\par
Assume that $\bar{G}$ is finite. Then
$J$ is infinite because so is $G$.
By Harbourne \cite{HPAMS} Proposition 1.3,
$X$ has only finitely many 'exceptional curves' (see \cite{HPAMS} page 409,
last paragraph for its definition)
and our $J$ stabilizes each irreducible component of these exceptional curves.
Inductively we can prove that there is a $J$-equivariant
smooth blowdown $X \rightarrow Y$ onto a relatively minimal rational surface $Y$
(with $K_Y^2 = 8, 9$).
Now the assertion (1) follows from the Riemann-Roch theorem.
\par
Assume $\bar{G}$ is infinite.
Then Theorem \ref{Thnull} follows from Theorem \ref{Thnulla}.
\par \vskip 0.5pc
Let $\sigma : \Aut(X) \rightarrow \Aut(\NS(X))$ be
the natural homomorphism.
\par \vskip 0.5pc
\begin{proposition}\label{seq}
Let $X$ be a smooth projective surface and $G \le \Aut(X)$
a subgroup of positive entropy.
Suppose that there is a sequence
$$H \triangleleft A \trianglelefteq G$$
such that $\sigma(H)$ (resp. $\sigma(A)$) is a finite (resp. infinte)
subgroup of $\Aut(\NS(X))$, and $A/H$ is abelian.
\par
Then there is a $B \le G$ fitting the sequence
below and satisfying the four conditions below
$$H \triangleleft A \trianglelefteq B \trianglelefteq G.$$
\begin{itemize}
\item[(1)]
$|G / B| \le 2$.
\item[(2)]
$B = \langle h_m \rangle \ltimes T$,
with $h_m$ positive entropy and $\sigma(T)$ a finite group.
\item[(3)]
$A = \langle h_m^a \rangle \ltimes T_1$,
with $a \ge 1$ and $T_1 \trianglelefteq T$.
\item[(4)]
$H \trianglelefteq T_1$ and
$A/H = \langle \bar{h}_m^a \rangle \times (T/H)$,
with $T/H$ abelian.
\end{itemize}
\end{proposition}
\begin{proof}
We claim that $A$ is of positive entropy.
Suppose the contrary that $A$ is of null entropy.
By Lemma \ref{Og2} and the proof of Theorem \ref{Thnulla}
($q(X) = 0$ was used only in the calculation of the Riemann Roch theorem),
there is a nef Cartier divisor $V$ such that $V^2 = 0$
and $a^*V \equiv V$ for all $a \in A$.
Now for $g \in G$ and $a \in A$, write
$g a g^{-1} = a' \in A$.
Then $(g a g^{-1})^*V = (a')^*V \equiv V$
and $a^*(g^*V) \equiv g^*V$ for all $a \in A$.
The uniqueness of the ray $\R_{>0} [V]$ in Lemma \ref{Og2},
implies $g^*V \equiv \alpha_g V$ for some $\alpha_g > 0$.
By the assumption on $G$,
some $g \in G$ is of positive entropy, whence
$L(g)^{\pm} = V$, a Cartier divisor (Lemma \ref{Lcor}). This contradicts
Proposition \ref{cyclo}.
Thus the claim is true.
\par
By Proposition \ref{nefstab3}, we have $A = \Stab_A(P)$
for every $P = L(g)^{\pm}$ with $g \in A$ being positive entropy.
Set $B := \Stab_G(P)$.
By Proposition \ref{nefstab2}, we have $|G : B| \le 2$
(hence $B \trianglelefteq G$).
Applying Theorem \ref{nefstab1} to $B$ and $A$,
we get (2) and (3).
For (4), we have $H \le T_1$ because $\sigma(H)$ is finite and
hence $H$ is of null entropy.
The proposition is proved.
\end{proof}
As a consequence of Proposition \ref{seq}, we have the result below:
\par \vskip 0.5pc
\begin{theorem}\label{Thz2a}
Let $X$ be a smooth projective surface and $G \le \Aut(X)$ a subgroup of positive
entropy.
Suppose that there is an $H \triangleleft G$ such that $G/H = \Z \times \Z$
and $\Imm (H \rightarrow \Aut(\NS(X)))$ is finite. Then we have:
\begin{itemize}
\item[(1)]
$G = \langle h_m \rangle \ltimes T$ and $H \trianglelefteq T$.
\item[(2)]
$\Imm (T \rightarrow \Aut(\NS(X)))$ is finite.
\item[(3)]
$G/H = \langle \bar{g}_1 \rangle \times \langle \bar{g}_2 \rangle$,
with $g_1 = h_m$ of positive entropy, and $g_2 \in T$ of null entropy and infinte order.
\end{itemize}
\end{theorem}
For groups which are not necessarily of null entropy,
we have Theorems \ref{Thsol} and \ref{Thcent} below, which are more general than
Theorem \ref{Thz2} in the Introduction.
\par
A group $F$ is {\it almost infinite cyclic} if
$F \cong \Z \ltimes F_1$ (semi-direct product) with $F_1$ finite (and normal in $F$).
\par \vskip 0.5pc
\begin{theorem}\label{Thsol}
Let $X$ be a smooth projective surface and $G \le \Aut(X)$ an infinite
subgroup. Assume the following two conditions:
\begin{itemize}
\item[(1)] There is an $H \trianglelefteq G$ such that
$G/H$ is soluble and the group $\Imm (H \rightarrow \Aut(\NS(X)))$
is finite; and
\item[(2)] No index $\le 2$ subgroup of $G$ is almost infinite cyclic.
\end{itemize}
Then $G$ contains a subgroup $S$ of null entropy and infinite order.
\newline
In particular, when $X$ is rational, there is an $S$-equivariant
smooth blowdown $X \rightarrow Y$ such that $Y$ has an anti-pluricanonical curve.
\end{theorem}
\begin{proof}
We may assume that $G$ is of positive entropy (see Theorem \ref{Thnull}).
Consider the derived series of $G$:
$$G = G^{(0)} \trianglerighteq G^{(1)} \trianglerighteq G^{(2)} \trianglerighteq \dots$$
with $G^{(n)} = [G^{(n-1)}, G^{(n-1)}]$
the commutator subgroup of $G^{(n-1)}$.
By the assumption, we have $G^{(n)} \le H$ for some $n \ge 0$,
whence $\sigma: \Aut(X) \rightarrow \Aut(\NS(X))$
maps $G^{(n)}$ to a finite group. Note that $\sigma(G)$ is infinte
because $G$ is of positive entropy by the additional assumption.
Thus there is a sequence below for some $r \ge 0$:
$$G^{(r+1)} \triangleleft G^{(r)} \trianglelefteq G$$
where $\sigma(G^{(r+1)})$ is finite, $\sigma(G^{(r)})$
is infinite, and $G^{(r)}/G^{(r+1)}$ is abelian.
\par
Applying Proposition \ref{seq} to $H = G^{(r+1)}$ and
$A = G^{(r)}$, we have $B = \Stab_G(P) = \langle h_m \rangle \ltimes T$
and $|G /B| \le 2$, as described there.
If $T$ is infinite, then let it be $S$ and we are done.
If $T$ is finite, then the index $\le 2$ subgroup $B$ of
$G$ is almost infinite cyclic, contradicting the assumption of the theorem.
This proves Theorem \ref{Thsol}.
\end{proof}
\par
For a group $P$, we consider the upper central series:
$$(1) = Z_0(P) \trianglelefteq Z_1(P) \trianglelefteq Z_2(P)
\trianglelefteq \cdots$$
with $Z_n(P)/Z_{n-1}(P)$ equal to
the centre $Z(P/Z_{n-1}(P))$.
\par \vskip 0.5pc
\begin{theorem}\label{Thcent}
Let $X$ be a smooth projective surface and $G \le \Aut(X)$ an infinite
subgroup. Assume the following two conditions:
\begin{itemize}
\item[(1)] There is an $H \triangleleft G$ such that the group
$\cup_{n=1}^{\infty} Z_n(G/H)$ is infinite
and the group $\Imm (H \rightarrow \Aut(\NS(X)))$
is finite; and
\item[(2)]
No index $\le 2$ subgroup of $G$ is almost infinite cyclic.
\end{itemize}
Then $G$ contains a subgroup $S$ of null entropy and infinite order.
\newline
In particular, when $X$ is rational, there is an $S$-equivariant
smooth blowdown $X \rightarrow Y$ such that $Y$ has an anti-pluricanonical curve.
\end{theorem}
\begin{proof}
We may assume that $G$ is of positive entropy (see Theorem \ref{Thnull}).
Consider the upper central series of $\bar{G} := G/H$:
$$(1) = Z_0(\bar{G}) \trianglelefteq Z_1(\bar{G}) \trianglelefteq Z_2(\bar{G})
\trianglelefteq \dots$$
with $Z_n(\bar{G})/Z_{n-1}(\bar{G}) = Z(\bar{G}/Z_{n-1}(\bar{G}))$,
the centre of $\bar{G}/Z_{n-1}(\bar{G})$.
Write $Z_n(\bar{G}) = Z_n/H$ for some $Z_n \trianglelefteq G$
with $Z_0 = H$.
By the assumption, $(\cup_{n \ge 1} Z_n)/H = \cup_{n \ge 1} Z_n(\bar{G})$ is
infinite. So $\cup_{n \ge 1} Z_n$ is infinite.
If $\sigma : \Aut(X) \rightarrow \Aut(\NS(X))$
maps every $Z_n$ ($n \ge 0$) to a finite group, then
the group $\cup_{n \ge 1} Z_n$ has null entropy and let it
be $S$, so we are done.
\par
Therefore, we assume that there is a sequence below for some $r \ge 1$:
$$Z_{r-1} \triangleleft Z_{r} \trianglelefteq G$$
where $\sigma(Z_{r-1})$ is finite, $\sigma(Z_r)$
is infinite, and $Z_r/Z_{r-1} \cong Z_r(\bar{G})/Z_{r-1}(\bar{G})
= Z(\bar{G}/Z_{r-1}(\bar{G}))$ is abelian.
The rest of the proof is identical to that of Theorem \ref{Thsol}.
This completes the proof of Theorem \ref{Thcent}.
\end{proof}
\begin{setup} {\bf Proof of Theorem \ref{Thz2}.}
\end{setup}
\par
This theorem will follow from Theorem \ref{Thsol} or \ref{Thcent},
but we give a direct proof,
so that we can see the dynamical structure of $G$.
If $G$ is of null entropy, we let it be $S$ and the theorem is true.
Assume that $G$ is of positive entropy.
We apply Proposition \ref{seq}
and use the notation there. Since $|G:A| = \infty$, we have
$|B : A| = \infty$, so $|T : T_1| = \infty$.
Thus $T$ is infinite and of null entropy; let it be $S$
and Theorem \ref{Thz2} is proved.
\par \vskip 0.5pc
\begin{setup} {\bf Proof of Theorem \ref{Thtits}.}
\end{setup}
\par
Set $H := G \cap \Ker(\sigma)$. By the original Tits alternative
theorem in \cite{Ti} Theorem 1, applied to $G/H \le \GL(\rho, \C)$
with $\rho = \rho(X)$, either $G/H \ge \Z * \Z = \langle \bar{g}_1
\rangle * \langle \bar{g}_2 \rangle$, the non-abelian free group of
rank $2$ (so $G \ge \Z * \Z = \langle g_1 \rangle * \langle g_2
\rangle$ and $G$ satisfies the assertion (1) of the theorem) , or
$G/H$ is virtually soluble: there is a $G_1/H \le G/H$ such that $|G
: G_1|$ is finite and $G_1/H$ is soluble.
\par
Let $\{g_iG_1\}$ be the (finite) set of all left cosets of $G_1$ in $G$.
Define $G_2 := \cap g_i G_1 g_i^{-1}$ ($\ge H$). Clearly, $G_2$ is normal in $G$.
Also $|G/G_2|$ is finite because all $|G : g_iG_1g_i^{-1}|$ are finite.
This can be proved inductively by the injectivity of the following map between the sets of
cosets:
$\{g (J_1 \cap J_2) \, | \, g \in G\} \rightarrow
\{g J_1 \, | \, g \in G\} \times \{g J_2 \, | \, g \in G\}$, \,
$x (J_1 \cap J_2) \mapsto (x J_1, \, x J_2)$; here $J_r \le G$.
\par
Since $G_1/H$ is soluble so is its subgroup $G_2/H$.
Since $G$ is of positive entropy, so is $G_2$ by the finiteness of $|G / G_2|$.
By the proof of Theorem \ref{Thsol} applied to $G_2$ and $H$,
there is a $B_2 = \Stab_{G_2}(P) = \langle h_2 \rangle \ltimes T_2$
such that $|G_2 / B_2| \le 2$.
\par
Set $B = \Stab_G(P)$.
If $B_2$ is normal in $G$, then the theorem (with the $B$ just defined)
follows from Proposition \ref{nefstab2} and Theorem \ref{nefstab1}.
\par
Suppose $B_2$ is not normal in $G$.
Then $|G_2/B_2| = 2$ and let $\tau \in G_2 \setminus B_2$.
We apply Proposition \ref{nefstab2} to $B_2 \trianglelefteq F := G_2 \le \Aut(X)$
and use its proof.
For any $g \in G$, \, $g h_2 g^{-1}$ is in $G_2$ and of positive entropy,
whence it equals $h_2^r t$ for some $r \in \{\pm 1\}$ and $t \in T_2$.
Also, if $r = 1$, then $g \in B$; if
$r = -1$ then $(\tau g) h_2 (\tau g)^{-1} = \tau(h_2^{-1} t) \tau^{-1} =
h_2 t_1 (\tau t \tau^{-1}) =: h_2 t_2$ with $t_1$ and $t_2$ in $T_2$,
so that $\tau g \in B$; also $\tau^2 \in T_2 \le B_2 \le B$.
Thus $G/B = \langle \bar{\tau} \rangle \cong \Z/(2)$.
Theorem \ref{Thtits} is proved (see Theorem \ref{nefstab1}).
\section{Positive entropy;  proof of Theorem \ref{Thpos-a}}
In this section, we consider surface automorphisms $g$ of positive entropy.
\par
{\it We prove} {\bf Theorem \ref{Thpos}} {\it below, which
will imply} {\bf Theorem \ref{Thpos-a}} {\it in the Introduction.}
We also show that the two sets $\Null(L(g)^{\pm})$
are identical (and finite)
and equal to $\Stab(g)$; see Lemma \ref{L3}.
\par \vskip 0.5pc
\begin{theorem}\label{Thpos}
Let $X$ be a smooth projective rational surface with
an automorphism $g$ of positive entropy.
Assume the following two conditions:
\begin{itemize}
\item[(1)] The pair $(X, \langle g \rangle)$ is minimal; and
\item[(2)] Either the set $\Stab(g)$ of $g$-periodic curves contains a curve
of arithmetic genus $\ge 1$,
or $X$ has an anti-pluricanonical curve.
\end{itemize}
Then we have:
\begin{itemize}
\item[(1)]
There is a unique nonzero effective $\Q$-divisor $\Delta$ with
$\Supp(\Delta) \subseteq \Stab(g)$ such that
$K_X + \Delta \equiv 0$. If $d$ is the smallest
positive integer such that $d\Delta$ is a Cartier integral divisor, then
$d(K_X + \Delta) \sim 0$.
Further, the $\gcd$ of coefficients of $d\Delta$ is equal to $1$.
\item[(2)] $\Stab(g)$ is a union of $\Supp(\Delta)$ and possibly
a few $(-2)$-curves away from $\Supp(\Delta)$.
\item[(3)]
One has $\kappa(X, -K_X) = 0$, so every anti-pluricanonical curve
is of the form $s(d\Delta)$ for some positive integer $s$.
\item[(4)]
Suppose that $d \ge 2$. Then $\Stab(g)$ is
a disjoint union of rational trees;
we have $C^2 \le -2$ for all $C \in \Stab(g)$ and
$C_1^2 \le -3$ for at least one $C_1 \in \Stab(g)$.
\end{itemize}
\end{theorem}
\par
We first prove the finiteness of
the sets $\Null(L^{\pm})$; the finiteness of
$\Null(L^+)$ $\cap \ \Null(L^-)$ is known and easy;
see section 2 for the notation.
\par \vskip 0.5pc
\begin{theorem}\label{infmany}
Let $X$ be a smooth projective surface with an automorphism $g$
of positive entropy.
Set $L^{\pm} = L(g)^{\pm}$.
\par
Then for each of $M = L^+$ and $M = L^-$ the set
$\Null(M)$ is either empty or a finite set.
In particular, $\Null(L^+) = \Null(L^-) = \Stab(g)$.
\end{theorem}
\begin{proof}
Clearly, there is a $\langle g \rangle$-equivariant
smooth blowdown $\gamma: X \rightarrow Y$ (with $E_{\gamma}$
the exceptional divisor) such that
$(Y, \langle g \rangle)$ is minimal.
Also, $E_{\gamma} \subseteq \Stab(g)$.
So $L^{\pm} . E_1 = 0$ for every $E_1 \le E_{\sigma}$ by Lemma \ref{L3}.
Thus $L^{\pm} = \gamma^*L_Y^{\pm}$ for some nef divisor
$L_Y^{\pm}$ on $Y$. The calculation (with $\lambda = \lambda(g)$)
$$\gamma^*(\lambda L_Y^{\pm}) = \lambda L^{\pm} \equiv g^* L^{\pm} = g^* \gamma^*L_Y^{\pm}
= \gamma^* (g^* L_Y^{\pm})$$
implies $g^*L_Y^{\pm} = \lambda_Y^{\pm}$ because
$\gamma^* : (\Pic(Y)) \otimes_{\Z} \R \rightarrow (\Pic(Y)) \otimes_{\Z} \R$
is injective. So $g$ is also of positive entropy on $Y$.
By the projection formula, we have:
$$\Null(L^{\pm}) = \gamma^{-1}(\Null(L^{\pm}_Y)) \cup \Supp(E_{\gamma}).$$
Thus we may assume that $(X, \langle g \rangle)$ is already minimal.
\par
Assume the contrary that $\Null(M)$ contains
infinitely many curves $C_i$. Since the Picard number $\rho(X)$ is finite,
there are positive integers $r$, $s$, $a_i$, $b_j$ such that:
$$D:= \sum_{i=1}^r a_i C_i \equiv \sum_{j=r+1}^{s} b_j C_j.$$
This display shows that $D$ is nef and is not numerically
trivial by the existence of an ample divisor $H$ on $X$.
Also $M . D = 0$. So $D \equiv a M$ for some $a > 0$ (Lemma \ref{nef}).
The Cartier-ness of $D$ contradicts Proposition \ref{cyclo}.
\end{proof}
We study $\Stab(g)$ according to the arithmetic genera of its members.
\par \vskip 0.5pc
\begin{lemma}\label{g=1}
Let $X$ be a smooth projective surface with an automorphism $g$
of positive entropy.
Let $C$ be a $g$-periodic curve.
Then we have:
\begin{itemize}
\item
[(1)]
The arithmetic genus $p_a(C) \le 1$.
\item
[(2)] If $p_a(C) = 1$ then $X$ is rational.
\item
[(3)] Suppose that $(X, \langle g \rangle)$ is minimal
and $p_a(C) = 1$. Then $C$ is an anti $1$-canonical curve.
$\Stab(g)$ is a union of $C$ and possibly a few $(-2)$-curves away from $C$.
\end{itemize}
\end{lemma}
\begin{proof}
By the assumption, we may fix some $s > 0$ such that
$g^s(C) = C$.
If $X$ is irrational, then $C$ is a smooth rational curve; see e.g.
Kawaguchi \cite{Kawaguchi} Proposition 3.1. In particular, (2) is true.
\par
We may assume that $X$ is rational.
We follow Diller, Jackson and Sommese \cite{DJS}, Theorem
3.6.
Considering the cohomology exact sequence associated to the exact
sequence below:
$$0 \rightarrow \OO_X(K_X) \rightarrow \OO_X(K_X + C)
\rightarrow \OO_C(K_C) \rightarrow 0,$$
we get $p_a(C) \le h^0(X, K_X + C)$ because $h^1(X, K_X) = q(X) = 0$.
If $|K_X + C| = \emptyset$ or
$K_X + C \sim 0$, (1) and (3) are true (see below).
\par
Suppose that
$0 < D \in |K_X + C|$. Then
$L . D = L . K_X + L . C = 0$ by Lemma \ref{L2}
and hence the intersection matrix of irreducible components
of $D$ is negative definite by the Hodge index theorem or Lemma \ref{L3}.
Thus $1 = h^0(X, D) \ge p_a(C)$ and also $g^s(D) = D$ as sets,
since $(g^s)^*D \sim D$ .
This proves the assertion (1).
\par
For (3), suppose $p_a(C) = 1$, $(X, \langle g \rangle)$ is minimal
and the above $D > 0$.
Note that $C$ is not a component of $D$, for otherwise $K_X \sim D - C \ge 0$,
absurd (by (2)).
Now $D^2 < 0$ and hence $0 > D . D_1 = (K_X + C ) . D_1 \ge
K_X . D_1$ for some component $D_1$ of $D$. Thus $D_1$ is a $(-1)$-curve.
Note that this $D_1 \in \Stab(g)$ because $g^s(D) = D$.
This contradicts the minimality of $(X, \langle g \rangle)$ by Lemma \ref{L3}.
Therefore, $K_X + C \sim 0$.
\par
To prove the last part of (3), let $C_1$ ($\ne C$) be in $\Stab(g)$.
Then $C_1^2 < 0$ by Lemma \ref{L3} and $0 = C_1 . (K_X + C) \ge C_1 . K_X$.
Thus $C_1$ is a $(-2)$-curve by the minimality of $(X, \langle g \rangle)$.
This proves the lemma.
\end{proof}
\begin{lemma}\label{finite}
Let $X$ be a smooth projective surface with an automorphism $g$
of positive entropy.
Suppose that $\Stab(g)$ is non-empty and consists of smooth rational curves.
Then we have:
\begin{itemize}
\item
[(1)] Every curve $C$ in $\Stab(g)$ is a smooth rational curve
with $C^2 = -n$ for some $n \ge 1$.
\item
[(2)] Suppose that $(X, \langle g \rangle)$ is minimal. Then in $(1)$ we have $n \ge 2$
and $K_X . C = n-2 \ge 0$. There is a unique effective $\Q$-divisor
$\Delta$ with $\Supp(\Delta) \subseteq \Stab(g)$
such that $(K_X + \Delta) . C_i = 0$ for every $C_i$ in $\Stab(g)$.
The set $\Stab(g)$ is a union of $\Supp(\Delta)$ and possibly
a few $(-2)$-curves away from $\Supp(\Delta)$.
\end{itemize}
\end{lemma}
\begin{proof}
(1) and the first part of (2) follow from Lemma \ref{L3}
while the second of (2) is from solving linear equations and Zariski's lemma
(Kollar and Mori \cite{KM} Lemma 3.41).
Take $C$ in $\Stab(g) \setminus \Supp(\Delta)$. Then
$C$ a $(-n)$-curve for some $n \ge 2$, and
we have $0 = C . (K_X + \Delta) \ge C . K_X = n-2 \ge 0$;
so the last of (2) follows. The lemma is proved.
\end{proof}
\begin{setup} {\bf Proof of Theorem \ref{Thpos}.}
\end{setup}
\par
If $\Stab(g)$ contains a curve $C$ with $p_a(C) \ge 1$,
then the theorem follows from Lemma \ref{g=1}.
\par
So we may assume that $X$ has an anti-pluricanonical curve
and $\Stab(g)$ (if not empty) consists of smooth rational curves $C_i$
(with $C_i^2 \le -2$ by Lemma \ref{finite}).
So $-d'K_X \sim D'$ for some $d' \ge 1$ and
effective (integral) Cartier divisor $D'$ ($\ne 0$, for $X$ being rational).
By Lemmas \ref{L2} - \ref{L3}, we have $D'. L(g)^{\pm} = -d'K_X . L(g)^{\pm} = 0$
and hence $\Supp(D') \subseteq \Stab(g)$.
Since $K_X + (D'/d') \sim_{\Q} 0$
and hence $(K_X + (D'/d')) . C_i = 0$ for every $C_i \in \Stab(g)$,
we have $\Delta = D'/d'$ by the uniqueness in Lemma \ref{finite}.
Also (2) follows.
\par
Write $\Delta = D/d$ with positive integer $d$ and effective Cartier
integral divisor $D$ such that the $\gcd$ (called $r$)
of coefficients of $D$ is coprime to $d$.
We claim that $r = 1$.
Indeed, take a $(-1)$-curve $E$ (noting that $K_X^2 < 0$ by Lemma \ref{L2}).
We have $1 = -K_X . E = (D . E)/d$.
Clearly, $r \, | \, D . E$. So $r \, | \, d$. Thus $r = 1$.
\par
Note that $dK_X + D = d(K_X + \Delta) \equiv 0$.
So $n(dK_X + D) \sim 0$ for some positive integer $n$,
since $q(X) = 0$. Thus $dK_X + D \sim 0$ because $X$ is
rational and hence $\Pic(X)$ is torsion free.
\par
Since $\Supp(\Delta) \subseteq \Stab(g)$ is negative definite, we have
$\kappa(X, -K_X) = 0$.
This proves (1). Now (3) follows.
\par
For (4), we assert that $C_{i_0}^2 \le -3$ for at least one $i_0$.
Otherwise, all $C_i$ are $(-2)$-curves and hence $K_X . \Delta = 0$.
This implies $0 = \Delta . (K_X + \Delta) = \Delta^2$,
whence $\Delta = 0$ by the negativity of $\Supp(\Delta) \subseteq \Stab(g)$.
So $K_X \equiv 0$, absurd (for $X$ being rational).
The assertion is proved.
\par
Suppose that a reduced connected component $D$ of $\Stab(g)$ is either of
non-simple-normal-crossing or contains a loop. Then
$(K_X + D) . D \ge 0$ (i.e., $p_a(D) \ge 1$); see e.g. CCZ \cite{CCZ} Lemma 2.2.
By the Riemann-Roch theorem, one has $K_X + D \sim G \ge 0$.
Cancelling the common components of $D$ and $G$, one has
$K_X + D' \sim G'$. If $G' = 0$ then one has a contradiction to that $d \ge 2$.
Thus $G' \ne 0$. By Lemma \ref{L2} - \ref{L3}, one has
$L . G' = L . (K_X + D') = 0$ and that $\Supp(G') \subseteq \Stab(g)$
is negative definite. Hence $(G')^2 < 0$.
So $0 > G_1 . G' = G_1 . (K_X + D') \ge G_1 . K_X$ for some component $G_1$
of $G'$. Thus $G_1$ is a $(-1)$-curve in $\Null(L)$,
contradicting the minimality of $(X, \langle g \rangle)$
(Lemma \ref{L3}). So (4) is true.
Theorem \ref{Thpos} is proved.
\par \vskip 0.5pc
As an effective upper bound for the $d$ in Theorem \ref{Thpos}, we have:
\par \vskip 0.5pc
\begin{remark}\label{Theff}
Let $X$ be a smooth projective rational surface with
an automorphism $g$ of positive entropy. Assume
the following three conditions:
\begin{itemize}
\item[(1)] The pair $(X, \langle g \rangle)$ is minimal;
\item[(2)] Either the set $\Stab(g)$ of $g$-periodic curves
contains a curve of arithmetic genus $\ge 1$,
or $X$ has an anti-pluricanonical curve; and
\item[(3)] Either $\Stab(g)$ is contractible to quotient singularities, or
$\Stab(g)$ is not a disjoint union of rational trees.
\end{itemize}
Then $X$ has an anti $s$-canonical curve
for some $1 \le s \le 21$.
\end{remark}
\par
Indeed,
As in Theorem \ref{Thpos}, $d(K_X + \Delta) \sim 0$.
We may assume that $d \ge 2$. Let $X \to \bar{X}$ be
the contraction of $\Delta$ to quotient singularities.
Then $\bar{X}$ is a rational log Enriques surface of index $d$
in the sense of \cite{Z1} Definitions 1.1 and 1.4.
Then by \cite{Z1} and \cite{Bl}, we have $d \le 21$.
\par \vskip 0.5pc
\begin{question}\label{ratvsirrat3}
Does there exist an upper bound of $s$ as in Remark $\ref{Theff}$ without
assuming the condition $(3) ? \ $
\end{question}
\end{large}
%
%

\end{document}